\documentclass{amsart}
\usepackage{amsfonts,amssymb,amsmath,amsthm}
\usepackage{url}
\usepackage{enumerate}

\newtheorem{theorem}{Theorem}[section]
\newtheorem{lemma}[theorem]{Lemma}
\newtheorem{prop}[theorem]{Proposition}
\newtheorem{corollary}[theorem]{Corollary}

\theoremstyle{definition}
\newtheorem{definition}[theorem]{Definition}
\newtheorem{example}[theorem]{Example}

\newtheorem*{notation}{Notation}

\renewcommand{\d}{\mathrm{d}}
\newcommand{\TD}{\mathrm{TD}}
\newcommand{\STD}{\mathrm{STD}}

\author{Robert F.\ Bailey}
\address{
School of Science \& Environment (Mathematics)\\
Grenfell Campus, Memorial University of Newfoundland\\
Corner Brook, NL~~A2H 6P9
}
\email{rbailey@grenfell.mun.ca}
\thanks{Supported by an NSERC Discovery Grant and a Memorial University of Newfoundland startup grant.}


\begin{document}

\title[Metric dimension of incidence graphs]{On the metric dimension of incidence graphs}

\begin{abstract}
A {\em resolving set} for a graph $\Gamma$ is a collection of vertices $S$, chosen so that for each vertex $v$, the list of distances from $v$ to the members of $S$ uniquely specifies $v$. The {\em \mbox{metric} dimension} $\mu(\Gamma)$ is the smallest size of a resolving set for $\Gamma$.  We consider the metric dimension of two families of incidence graphs: incidence graphs of symmetric designs, and incidence graphs of symmetric transversal designs (i.e.\ symmetric nets).  These graphs are the bipartite distance-regular graphs of diameter~$3$, and the bipartite, antipodal distance-regular graphs of diameter~$4$, respectively.  In each case, we use the probabilistic method in the manner used by Babai to obtain bounds on the metric dimension of strongly regular graphs, and are able to show that $\mu(\Gamma)=O(\sqrt{n}\log n)$ (where $n$ is the number of vertices).
\end{abstract}

\maketitle

\section{Introduction} \label{section:intro}

We consider finite, connected graphs with no loops or multiple edges.  Let $\Gamma$ denote a graph with vertex set $V$ and edge set $E$.  A {\em resolving set} for $\Gamma$ is a subset $S\subseteq V$ with the property that, for any $u\in V$, the list of distances from $u$ to each of the elements of $S$ uniquely identifies $u$; equivalently, for two distinct vertices $u,w\in V$, there exists $x\in S$ for which $\d(u,x)\neq\d(w,x)$.  The {\em metric dimension} of $\Gamma$ is the smallest size of a resolving set for $\Gamma$, and we denote this by $\mu(\Gamma)$.  These notions were introduced to graph theory in the 1970s by Slater~\cite{Slater75} and, independently, Harary and Melter~\cite{Harary76}; in more general metric spaces, the concept can be found in the literature much earlier (see~\cite{Blumenthal53}).  For further details, the reader is referred to the survey~\cite{bsmd}.

When studying metric dimension, distance-regular graphs are a natural class of graphs to consider.  A graph $\Gamma$ with diameter~$d$ is {\em distance-regular} if, for all $i$ with $0\leq i\leq d$ and any vertices $u,w$ with $\d(u,w)=i$, the number of neighbours of $w$ at distances $i-1$, $i$ and $i+1$ from $u$ depend only on the distance $i$, and not on the choices of $u$ and $w$.  These numbers are denoted by $c_i$, $a_i$ and $b_i$ respectively, and are known as the {\em parameters} of $\Gamma$.  It is easy to see that $c_0$, $b_d$ are undefined, $a_0=0$, $c_1=1$ and $c_i+a_i+b_i=k$ (where $k$ is the valency of $\Gamma$).  We put the parameters into an array, called the {\em intersection array} of $\Gamma$,
\[ \left\{ \begin{array}{cccccc}
\ast & 1   & c_2 & \cdots & c_{d-1} & c_d \\ 
 0   & a_1 & a_2 & \cdots & a_{d-1} & a_d \\
 k   & b_1 & b_2 & \cdots & b_{d-1} & \ast 
\end{array} \right\}. \]
Since the 2011 survey article by Cameron and the present author~\cite{bsmd}, which first proposed its systematic study, a number of papers have been written on the subject of the metric dimension of distance-regular graphs (and on the related problem of class dimension of association schemes), by the present author and others: see 
\cite{small,imprimitive,jk,grassmann,Bartoli,Beardon,FengWang,Gravier,attenuated,flat,dualpolar,GuoWangLi,fourfamilies,Heger}, for instance; earlier results may be found in \cite{Babai80,Babai81,Caceres07,Chvatal83,SeboTannier04}.  For background on distance-regular graphs in general, see the book of Brouwer, Cohen and Neumaier~\cite{BCN} or the survey by van Dam, Koolen and Tanaka~\cite{vanDam}.

A distance-regular graph $\Gamma$ with diameter~$d$ is {\em primitive} if, for $1\leq i\leq d$, the distance-$i$ graphs of $\Gamma$ are all connected; otherwise, we say it is {\em imprimitive}.  Imprimitive distance-regular graphs arise in one of two ways, provided that the valency is at least~$3$: they may be bipartite (whereby the distance-$2$ graph has two connected components, called the {\em halved graphs} of $\Gamma$), or {\em antipodal} (where the distance-$d$ graph is a disjoint union of cliques).  We note that both possibilities may occur in the same graph.  Imprimitive distance-regular graphs may be reduced to primitive ones by the operations of {\em halving} (for bipartite graphs) or {\em folding} (for antipodal graphs); see~\cite[{\S}4.2A]{BCN} for details.  If an imprimitive graph $\Gamma$ has diameter~$d\geq 3$, its halved or folded graphs have diameter $\lfloor d/2 \rfloor$.

The metric dimension of imprimitive distance-regular graphs was studied in detail in~\cite{imprimitive}, where it was shown that it can be bounded in terms of the metric dimension of the halved or folded graphs (see~\cite[{\S}2.1]{imprimitive}).  However, when the halved or folded graphs are either complete or complete multipartite, the results are unsatisfactory; this is especially true from the asymptotic perspective, as we obtain the trivial upper bound of $O(n)$ (where $n$ is the number of vertices).  In this paper, we consider bipartite distance-regular graphs of diameter~$3$, and distance-regular graphs of diameter~$4$ which are both bipartite and antipodal.  The former class is precisely equivalent to the incidence graphs of {\em symmetric designs}, which are well-understood objects (see~\cite{Ionin06}, for instance); the latter class is equivalent to the incidence graphs of {\em symmetric transversal designs}, or equivalently {\em symmetric nets}, about which the literature is more sporadic.

\subsection{Split resolving sets and semi-resolving sets}
In~\cite{imprimitive}, the present author introduced the following special type of resolving set for bipartite graphs.

\begin{definition} \label{defn:split-semi}
Let $\Gamma$ be a bipartite graph, whose vertex set has bipartition $X\cup Y$.  A {\em split resolving set} for $\Gamma$ is a subset of vertices $S=S_X\cup S_Y$, where $S_X\subseteq X$ and $S_Y\subseteq Y$, chosen so that any two vertices in $X$ are resolved by a vertex in $S_Y$, and any two vertices in $Y$ are resolved by a vertex in $S_X$.  We call $S_X$ a {\em semi-resolving set} for $Y$ and $S_Y$ a semi-resolving set for $X$.  We denote the smallest size of a split resolving set by $\mu^\ast(\Gamma)$.
\end{definition}

We note that a split resolving set is itself a resolving set: any vertex of $\Gamma$ will resolve a pair of vertices $(x,y)$ where $x\in X$ and $y\in Y$, given that the parities of the distances to $x$ and to $y$ will be different, so we only need consider resolving pairs of vertices in the same bipartite half.  Consequently, we have $\mu(\Gamma)\leq \mu^\ast(\Gamma)$.  We also note that complete bipartite graphs do not have split resolving sets.

If we regard a bipartite graph $\Gamma$ as an incidence graph, semi-resolving sets are of independent interest due to connections with other objects associated with incidence structures, such as blocking sets in finite geometries; see~\cite{imprimitive,Bartoli,Heger} for more details on this.

\section{Symmetric designs} \label{section:symmdes}

A {\em symmetric design} (or {\em square $2$-design}) with parameters $(v,k,\lambda)$ is a pair $\mathcal{D}=(X,\mathcal{B})$, where $X$ is a set of $v$ \emph{points}, and $\mathcal{B}$ is a family of $k$-subsets of $X$, called \emph{blocks}, such that any pair of distinct points are contained in exactly $\lambda$ blocks, and that any pair of distinct blocks intersect in exactly $\lambda$ points.  It follows that $|\mathcal{B}|=v$.  A symmetric design with $\lambda=1$ is a {\em projective plane}, while a symmetric design with $\lambda=2$ is known as a {\em biplane}.  The {\em incidence graph} $\Gamma_\mathcal{D}$ of a symmetric design $\mathcal{D}$ is the bipartite graph with vertex set $X\cup\mathcal{B}$, with the point $x\in X$ adjacent to the block $B\in\mathcal{B}$ if and only if $x\in B$.  It is straightforward to show that the incidence graph of a symmetric design is a bipartite distance-regular graph with diameter~$3$ and intersection array
\[ \left\{ \begin{array}{cccc}
  \ast & 1   & \lambda   & k \\
	0    & 0   & 0         & 0 \\
	k    & k-1 & k-\lambda & \ast
	\end{array} \right\}. \]
The converse is also true (see~\cite[{\S}1.6]{BCN}): any bipartite distance-regular graph of diameter~$3$ gives rise to a symmetric design.  The {\em dual} of a symmetric design is the design obtained from the incidence graph by reversing the roles of points and blocks; both $\mathcal{D}$ and its dual have the same parameters.

The {\em order} of a symmetric design is defined to be $q=k-\lambda$; the following result is well-known (see~\cite[Proposition 2.4.12]{Ionin06}, for instance) and gives restrictions on $v$ in terms of the order.

\begin{prop} \label{prop:order}
For any $(v,k,\lambda)$ symmetric design of order $q=k-\lambda \geq 2$, we have
\[ 4q-1 \leq v \leq q^2+q+1. \]
\end{prop}

The two extremes are achieved by Hadamard designs (where $v=4q-1)$ and projective planes (where $v=q^2+q+1)$.

The incidence graphs of symmetric designs are precisely the bipartite distance-regular graphs of diameter~$3$; the metric dimension of these graphs is considered in~\cite{imprimitive}.  However, the general results of~\cite{imprimitive} for bipartite distance-regular graphs are not very effective in the diameter~$3$ case, as the halved graphs are complete graphs, so an alternative approach was required.  First, in the case where $k=v-1$ (or, equivalently, where the order is $q=k-\lambda=1$), the incidence graph is $K_{v,v}-I$, i.e.\ a complete bipartite graph with a perfect matching removed, which has metric dimension $v-1$ (see~\cite[Corollary~2.7]{imprimitive}), so is linear in the number of vertices.  In the case where $k<v-1$, it was shown that if the design has a null polarity (an incidence-preserving bijection $\sigma$ between the points and blocks, where no point is incident with its image under $\sigma$), then $\mu(\Gamma)=O(\sqrt{n}\log n)$ (where $n$ is the number of vertices; see~\cite[Corollary~4.6]{imprimitive}).  Further, in~\cite{Heger}, H\'eger and Tak\'ats considered the special case of projective planes: they determined the exact size of a resolving set for any sufficiently large projective plane, and also the exact size of a split resolving set in the case of the Desarguesian plane $\mathrm{PG}(2,q)$.  Asymptotically, both values are $\Theta(\sqrt{n})$.

In this section, we will consider the metric dimension of incidence graphs of symmetric designs in general.  We shall do so by demonstrating the existence of a semi-resolving set of an appropriate size for the points of $\mathcal{D}$.  The following notation is useful.

\begin{notation}
Let $B(x)$ denote the set of blocks containing the point $x$.
\end{notation}

The next lemma provides a characterization of semi-resolving sets for the points of $\mathcal{D}$.  (We let $A\vartriangle B$ denote the symmetric difference of sets $A$ and $B$.)

\begin{lemma} \label{lemma:split}
For a symmetric design $\mathcal{D}=(X,\mathcal{B})$, a subset $S_\mathcal{B}\subseteq \mathcal{B}$ is a semi-resolving set for the points of $\mathcal{D}$ if and only if, for all pairs of distinct points $x,y\in X$, we have $S_\mathcal{B}\cap (B(x)\vartriangle B(y)) \neq \emptyset$.
\end{lemma}

\proof  We note that in $\Gamma_\mathcal{D}$, we have $\d(x,B)=1$ if and only if $x\in B$, and $\d(x,B)=3$ if and only if $x\not\in B$.  Thus for a block $B$ to resolve a pair of points $x,y$, we must have that it is incident with exactly one of the points $x,y$.  So, in order for $S_\mathcal{B}$ to be a semi-resolving set, it must contain a block incident with exactly one of $x,y$.  That is, $S_\mathcal{B}\cap (B(x)\vartriangle B(y)) \neq \emptyset$.
\endproof

Fortunately, the size of $B(x)\vartriangle B(y)$ is easy to calculate in terms of the parameters of $\mathcal{D}$.

\begin{lemma} \label{lemma:symmdiffsize}
Let $\mathcal{D}=(X,\mathcal{B})$ be a $(v,k,\lambda)$ symmetric design.
For any distinct points $x,y\in X$, we have $|B(x)\vartriangle B(y)|=2(k-\lambda)$.
\end{lemma}

\proof Clearly, $|B(x)|=|B(y)|=k$ and $|B(x)\cap B(y)|=\lambda$, and consequently $|B(x)\vartriangle B(y)|=|B(x)|+|B(y)|-2|B(x)\cap B(y)|=2(k-\lambda)$. \endproof

With these lemmas, we now present our main result of this section.

\begin{theorem} \label{thm:semiresolving}
Let $\mathcal{D}=(X,\mathcal{B})$ be a $(v,k,\lambda)$ symmetric design of order $q\geq 2$.  Then there exists a semi-resolving set for the points of $\mathcal{D}$ of size $\displaystyle \left\lceil \frac{v\log v}{k-\lambda} \right\rceil$.
\end{theorem}

\proof Let $m=2q=2(k-\lambda)$, so that by Lemma~\ref{lemma:symmdiffsize}, $m=|B(x)\vartriangle B(y)|$ for any distinct $x,y\in X$.

We use the probabilistic method (described in detail in Alon and Spencer~\cite{AlonSpencer08}).  
First, we let $s = \left\lceil v\log v/(k-\lambda) \right\rceil = \left\lceil 2v\log v/m \right\rceil$.  We note that $s\leq v$ if and only if $v\leq e^q$; since $q\geq 2$, Proposition~\ref{prop:order} implies that $v\leq q^2+q+1 <e^q$.

Suppose that $S_\mathcal{B}$ is a set of $s$ blocks chosen uniformly at random from $\mathcal{B}$, each set chosen with probability $1/ \binom{v}{s}$.  For any distinct $x,y\in X$, from Lemma~\ref{lemma:split} we know that $S_\mathcal{B}$ resolves the pair $\{x,y\}$ if and only if $S_\mathcal{B} \cap (B(x)\vartriangle B(y))$ is non-empty.  Let $A(x,y)$ denote the event that $S_\mathcal{B}$ {\em fails} to resolve $\{x,y\}$, and let $\mathrm{P}(x,y)=\mathrm{Pr}(A(x,y))$.  Since $|B(x)\vartriangle B(y)|=m$, we have that 
\[ \mathrm{P}(x,y) = \mathrm{Pr} \left( S_\mathcal{B}\cap (B(x)\vartriangle B(y)) = \emptyset \right) = \binom{v-m}{s} \biggm/ \binom{v}{s}. \]
Let $N$ denote the number of pairs $\{x,y\}$ such that $A(x,y)$ holds.  The expected value of $N$ is therefore
\[ \mathbb{E}(N) = \sum_{ \{x,y\}\subseteq X } \mathrm{P}(x,y) = \binom{v}{2}\binom{v-m}{s} \biggm/ \binom{v}{s}. \]
Clearly, $N=0$ if and only if $S_\mathcal{B}$ is a semi-resolving set, so if we can show that $\mathbb{E}(N) < 1$ there must be a semi-resolving set of size $s$. If $s>v-m$, then $P(x,y)$ is always $0$ and thus $\mathbb{E}(N)=0$, so we will assume that $s\leq v-m$.  In that case, we have $\mathbb{E}(N) < 1$ if and only if 
\[ \binom{v}{2}\binom{v-m}{s} < \binom{v}{s}. \]
Since $s=\left\lceil 2v\log v/m \right\rceil$, we have
\begin{align*}
      & \frac{2v\log v - v\log 2}{m} < s \\
 \iff & 2\log v - \log 2 < \frac{ms}{v} \\
 \iff & \frac{v^2}{2} < \exp(m/v)^s.
\end{align*}
From this, it follows that
\[ \binom{v}{2} < \frac{v^2}{2} < \exp(m/v)^s < \left( 1 + \frac{m}{v} + \frac{m^2}{v^2} \right)^s < \prod_{i=0}^{s-1} \left( 1+ \frac{m}{v-m-i} \right) = \binom{v}{s} \biggm/ \binom{v-m}{s}. \]
(Note that $e^t< 1+t+t^2$ for $0<t<1$; the final inequality can be proved by induction on $s$.)  Hence $\mathbb{E}(N) < 1$, and a semi-resolving set of size $s$ is guaranteed to exist.
\endproof

We remark that this method of proof closely follows that of Babai~\cite[Lemma~3.2]{Babai80}.  Theorem~\ref{thm:semiresolving} has the following immediate consequence.

\begin{corollary} \label{cor:symmdesresolving}
Let $\Gamma_{\mathcal{D}}$ denote the incidence graph of the $(v,k,\lambda)$ symmetric design $\mathcal{D}=(X,\mathcal{B})$, which has of order $q\geq 2$.  Then there exists a split resolving set for $\Gamma_{\mathcal{D}}$ of size 
$\displaystyle 2\cdot \left\lceil \frac{v\log v}{k-\lambda} \right\rceil$.
\end{corollary}


To analyse this result asymptotically, we appeal to Proposition~\ref{prop:order} which gives a lower bound on the order $q=k-\lambda$ of $\Omega(\sqrt{v})$; this yields the following result.

\begin{corollary} \label{cor:symmdesasymptotic}
Let $\Gamma$ be the incidence graph of a $(v,k,\lambda)$ symmetric design of order $q\geq 2$, and let $n=2v$ be the number of vertices.  Then $\mu(\Gamma)=O(\sqrt{n}\log n)$.
\end{corollary}

However, thanks to the lower bound on $v$ in Proposition~\ref{prop:order}, in the case where $v$ is linear in $q$ (such as for incidence graphs of Hadamard designs) we can obtain a tighter result.  When combined with the general observation that, for any graph $\Gamma$ on $n$ vertices, the metric dimension satisfies $\mu(\Gamma)=\Omega(\log n)$ (see~\cite[Proposition 3.6]{bsmd}), we have the following result.

\begin{corollary} \label{cor:strongerasymptotic}
Let $\Gamma$ be the incidence graph of a $(v,k,\lambda)$ symmetric design where $v$ is a linear function of $q=k-\lambda$, and let $n=2v$ be the number of vertices of $\Gamma$.  Then $\mu(\Gamma)=\Theta(\log n)$.
\end{corollary}

In~\cite{imprimitive}, the same bounds on $\mu(\Gamma)$ as found in Corollaries~\ref{cor:symmdesasymptotic} and~\ref{cor:strongerasymptotic} were obtained, but in the special case that the design has a null polarity (an incidence-preserving bijection $\sigma$ between points and blocks, where no point is incident to its image under $\sigma$).  These results remove that additional assumption, and thus answering an open question from~\cite{imprimitive} in the affirmative.  Corollary~\ref{cor:symmdesasymptotic} also extends to the class of bipartite distance-regular graphs of diameter~$3$ a result of Babai~\cite{Babai81}, which asserts that $\mu(\Gamma)=O(\sqrt{n}\log n)$ for primitive distance-regular graphs on $n$ vertices.


\section{Symmetric transversal designs} \label{section:std}

A {\em transversal design} $\TD_\lambda[k;g]$ is a triple $\mathcal{D}=(X,\mathcal{G},\mathcal{B})$, where $X$ is a set of $v=kg$ points, $\mathcal{G}$ is a partition of $X$ into $k$ sets of size $g$ called {\em point classes} (or, often in the literature, ``groups''), and $\mathcal{B}$ is a family of $k$-subsets of $X$ called blocks, chosen so that each block contains exactly one element of each point class, and so that any pair of points from distinct point classes lie in exactly $\lambda$ blocks.  Since a $\TD_\lambda[k;1]$ consists of a single block repeated $\lambda$ times, we assume that $g\geq 2$.  General background information on transversal designs can be found in Beth, Jungnickel and Lenz~\cite{BJL}. 

A {\em symmetric transversal design}, denoted $\STD_\lambda[k;g]$, is a $\TD_\lambda[k;g]$ whose dual (i.e.\ the incidence structure obtained by interchanging the roles of points and blocks) is also a $\TD_\lambda[k;g]$.  In such a design, it follows that $k=\lambda g$ and $v=|\mathcal{B}|=\lambda g^2$, so the parameters depend only on $\lambda$ and $g$.  The dual of a transversal design is often called a {\em net}; for this reason, symmetric transversal designs are often known as {\em symmetric nets}.  A $\TD_\lambda[k;g]$ is {\em resolvable} if its set of blocks $\mathcal{B}$ can be partitioned into {\em parallel classes} (or {\em block classes}), each of which is a partition of the point set $X$.  By definition, a symmetric transversal design is resolvable (the point classes give the parallel classes of the dual, and vice-versa).  Conversely, a resolvable transversal design with $k=\lambda g$ is necessarily symmetric (see~\cite[Theorem II.8.21]{BJL}).

The {\em incidence graph} $\Gamma_\mathcal{D}$ of a symmetric transversal design $\mathcal{D}$ is defined as in the previous section: we have a bipartite graph with vertex set $X\cup\mathcal{B}$, with $x\in X$ adjacent to $B\in\mathcal{B}$ if and only if $x\in B$.  The incidence graph of an $\STD_\lambda[k;g]$ is a bipartite, antipodal distance-regular graph with diameter~$4$ and intersection array
\[ \left\{ \begin{array}{ccccc}
  \ast      & 1           & \lambda      & \lambda g-1 & \lambda g \\
	0         & 0           & 0            & 0           & 0 \\
	\lambda g & \lambda g-1 & \lambda(g-1) & 1           & \ast
	\end{array} \right\}. \]
The converse is also true: any bipartite, antipodal distance-regular graph with diameter~$4$ is the incidence graph of a symmetric transversal design (see~\cite[{\S}1.7]{BCN}).  Their halved graphs are complete multipartite, while their folded graphs are complete bipartite; consequently, the results of~\cite{imprimitive} using halving or folding to obtain bounds on metric dimension do not yield particularly useful results here.  However, the methods of the previous section, using semi-resolving sets and split resolving sets, may be easily adapted to obtain better bounds on $\mu(\Gamma_\mathcal{D})$.  

Semi-resolving sets for the points of an $\STD_\lambda[k;g]$ may be characterized in exactly the same way as for symmetric designs.  As before, we let $B(x)$ denote the set of blocks containing the point $x$.

\begin{lemma} \label{lemma:split2}
For a symmetric transversal design $\mathcal{D}=(X,\mathcal{G},\mathcal{B})$, a subset $S_\mathcal{B}\subseteq \mathcal{B}$ is a semi-resolving set for the points of $\mathcal{D}$ if and only if, for all pairs of distinct points $x,y\in X$, we have $S_\mathcal{B}\cap (B(x)\vartriangle B(y)) \neq \emptyset$.
\end{lemma}

\proof Identical to the proof of Lemma~\ref{lemma:split}. \endproof

Unlike the case of symmetric designs, the size of $B(x)\vartriangle B(y)$ is not constant for all pairs of distinct points $x,y\in X$.  However, there are only two possibilities, depending on whether $x$ and $y$ lie in the same point class or not.

\begin{lemma} \label{lemma:symmdiffsize2}
Let $\mathcal{D}=(X,\mathcal{G},\mathcal{B})$ be an $\STD_\lambda[k;g]$.  For distinct points $x,y\in X$, we have 
\[ |B(x)\vartriangle B(y)| = \left\{ \begin{array}{ll}
                                  2k = 2\lambda g              & \textnormal{if $x,y$ lie in the same point class,} \\
																	2(k-\lambda) = 2\lambda(g-1) & \textnormal{if $x,y$ lie in distinct point classes.}
                                  \end{array} \right. \]
\end{lemma}

\proof If $x,y$ lie in the same point class, we have $B(x)\cap B(y)=\emptyset$, while if $x,y$ lie in distinct point classes we have $|B(x)\cap B(y)|=\lambda$; the result follows from this. \endproof

The main result of this section bears a strong similarity to Theorem~\ref{thm:semiresolving} for symmetric designs.

\begin{theorem} \label{thm:semiresolving2}
Let $\mathcal{D}=(X,\mathcal{G},\mathcal{B})$ be an $\STD_\lambda[k;g]$, so $v=\lambda g^2$ and $k=\lambda g$, where $\lambda\geq 1$, $g\geq 2$ and $(\lambda,g)\not\in \{(1,2),(1,3),(2,2)\}$.  Then there exists a semi-resolving set for the points of $\mathcal{D}$ of size $\displaystyle \left\lceil \frac{v\log v}{k-\lambda} \right\rceil$.
\end{theorem}

\proof As with the proof of Theorem~\ref{thm:semiresolving}, we will use the probabilistic method.  Once again, our aim is to show that the expected number $N$ of pairs of points $\{x,y\}$ that fail to be resolved by a set of $s$ blocks $S_\mathcal{B}$, chosen uniformly at random from $\mathcal{B}$, is less than $1$.  Let $s = \left\lceil v\log v/(k-\lambda) \right\rceil$; the restrictions $\lambda\geq 1$, $g\geq 2$ and $(\lambda,g)\not\in \{(1,2),(1,3),(2,2)\}$ ensure that $s\leq v$.

As before, we let $\mathrm{P}(x,y)$ denote the probability that $S_\mathcal{B}$ fails to resolve the pair $\{x,y\}$, i.e.\
\[ \mathrm{P}(x,y) = \mathrm{Pr} \left( S_\mathcal{B} \cap (B(x)\vartriangle B(y)) = \emptyset \right). \]
By Lemma~\ref{lemma:symmdiffsize2}, this probability will be
\[ \mathrm{P}(x,y) = \left\{ \begin{array}{ll}
                             {\displaystyle \binom{v-2k}{s} \biggm/ \binom{v}{s} }           & \textnormal{if $x,y$ lie in the same point class,} \\[3ex]
														 {\displaystyle \binom{v-2(k-\lambda)}{s} \biggm/ \binom{v}{s} } & \textnormal{if $x,y$ lie in distinct point classes.}
                             \end{array} \right. \]
Using the larger of these two possible values, we therefore obtain
\[ \mathbb{E}(N) = \sum_{ \{x,y\}\subseteq X } \mathrm{P}(x,y) \leq \binom{v}{2} \binom{v-2(k-\lambda)}{s} \biggm/ \binom{v}{s}. \]
By the same calculations as in the proof of Theorem~\ref{thm:semiresolving}, this yields $\mathbb{E}(N)<1$, so there must exist a semi-resolving set of size $s$.
\endproof

\begin{corollary} \label{cor:STDresolving}
Let $\mathcal{D}=(X,\mathcal{G},\mathcal{B})$ be an $\STD_\lambda[k;g]$, so $v=\lambda g^2$ and $k=\lambda g$, where $\lambda\geq 1$, $g\geq 2$ and $(\lambda,g)\not\in \{(1,2),(1,3),(2,2)\}$, and let $\Gamma_\mathcal{D}$ be its incidence graph.  Then there exists a split resolving set for $\Gamma_{\mathcal{D}}$ of size 
$\displaystyle 2\cdot\left\lceil \frac{v\log v}{k-\lambda} \right\rceil$.
\end{corollary}

We remark that the graphs in the three exceptional cases are the $8$-cycle, Pappus graph and $4$-cube, which have metric dimension $2$, $4$ and $4$ respectively (see~\cite{small}).  

In terms of the parameters $g$ and $\lambda$, Corollary~\ref{cor:STDresolving} gives a bound on the metric dimension of $\Gamma_\mathcal{D}$:
\[ \mu(\Gamma_\mathcal{D}) \leq 2\cdot\left\lceil \frac{g^2\log (\lambda g^2)}{g-1} \right\rceil, \]
which asymptotically behaves as $O(g (\log g + \log\lambda))$.  This leads to the following corollary.

\begin{corollary} \label{cor:STDasymptotic}
Let $\Gamma_\mathcal{D}$ be the incidence graph of an $\STD_\lambda[k;g]$, with $n=2\lambda g^2$ vertices.  Then $\mu(\Gamma_\mathcal{D})=O(\sqrt{n}\log n)$.
\end{corollary}


However, as with the incidence graphs of symmetric designs, we can obtain more precise results for specific families of graphs.  Once again, there are extreme cases which arise from projective geometries at one end, and Hadamard matrices at the other.  First, we consider the case where $\lambda=1$ and $g$ is allowed to vary.

\begin{example} \label{example:biaffine}
Let $g=q$ and $\lambda=1$.  An $\STD_1[q;q]$ is equivalent to an affine plane of order $q$ with a parallel class removed (sometimes called a {\em biaffine plane}), and thus exists whenever there exists a projective plane of order~$q$: see~\cite[Proposition II.7.19]{BJL}.  By Theorem~\ref{thm:semiresolving2}, there exists a semi-resolving set for a biaffine plane of order $q$ with size $\left\lceil 2q^2\log q /(q-1) \right\rceil$.  The incidence graph of a biaffine plane thus has metric dimension at most $2\cdot\left\lceil 2q^2\log q /(q-1) \right\rceil$.  However, a stronger result is possible: in~\cite[Propositions~3.2, 3.4]{Bartoli}, Bartoli {\em et al.}\ prove that if $\mathcal{D}$ is a biaffine plane of order~$q\geq 4$, then its incidence graph $\Gamma_\mathcal{D}$ satisfies
\[ 2q-2 \leq \mu(\Gamma_\mathcal{D}) \leq 3q-6. \]
Asymptotically, this shows that $\mu(\Gamma_\mathcal{D})=\Theta(\sqrt{n})$.  The upper bound is obtained by an explicit construction of a resolving set.  A tighter lower bound of $8q/3 -7$ is possible in the case of the Desarguesian biaffine plane $\mathrm{BG}(2,q)$ for $q\geq 7$ (see~\cite[Theorem~3.20]{Bartoli}), and with additional restrictions on~$q$ the lower bound may be improved further to $3q-9\sqrt{q}$ (see~\cite[Theorem~3.16]{Bartoli}), but neither of these affect the asymptotic result.
\end{example}

On the other hand, if $g$ is fixed and $\lambda$ is allowed to vary, the smallest interesting case is $g=2$.

\begin{example} \label{example:hadamardgraphs}
An $\STD_\lambda[2\lambda;2]$ is equivalent to a Hadamard matrix of order $2\lambda$ (and thus either $\lambda=1$ or $\lambda$ must be even); the corresponding incidence graph is called a {\em Hadamard graph} of order $2\lambda$ (see~\cite[{\S}1.8]{BCN}).  By Corollary~\ref{cor:STDresolving}, a Hadamard graph of order $2\lambda$ has a split resolving set of size $2\lceil 4\log(4\lambda) \rceil$.  Asymptotically, since a Hadamard graph $\Gamma$ has $n=8\lambda$ vertices, this---along with the general logarithmic lower bound on metric dimension---shows that $\mu(\Gamma)=\Theta(\log n)$.
\end{example}


\section{Conclusion} \label{section:conclusion}

In~\cite[Theorem 2.9]{AH06}, Alfuraidan and Hall gave a result which groups distance-regular graphs into various classes, in terms of whether the graph is primitive, bipartite, antipodal (or both bipartite and antipodal), and its diameter.  In~\cite{imprimitive}, these classes were labelled AH1--AH13: class AH1 consists of the primitive graphs of diameter at least~$2$ and valency at least~$3$; classes AH2--AH4 consist of cycles, complete graphs and complete multipartite graphs respectively; class AH5 consists of the graphs obtained by deleting a perfect matching from a complete bipartite graph.  The remaining classes consist of the remaining antipodal graphs; in the present paper, we are primarily concerned with the following two of them:
\begin{itemize}
  \item[AH6:] $\Gamma$ has diameter~$3$, is bipartite but not antipodal, its halved graphs are complete graphs, and $\Gamma$ is the incidence graph of a $(v,k,\lambda)$ symmetric design with block size $k<v-1$;
	\item[AH8:] $\Gamma$ has diameter~$4$, is both bipartite and antipodal, its halved graphs are complete multipartite, its folded graph is complete bipartite, and $\Gamma$ is the incidence graph of an $\STD_\lambda[k;g]$.
\end{itemize}
The results of Sections~\ref{section:symmdes} and~\ref{section:std} show that graphs in both of these classes have metric dimension $\mu(\Gamma)=O(\sqrt{n}\log n)$ (where $n$ is the number of vertices).  This extends Babai's results~\cite{Babai80,Babai81}, which give the same upper bound for primitive graphs (which form class AH1), and the results of~\cite{imprimitive}, which do so for graphs in classes AH11--AH13.\footnote{Note that there is a typographical error in the conclusion of~\cite{imprimitive}: it should refer to classes AH10--AH13, rather than AH9--AH12.} 

The class AH7, which consists of the non-bipartite antipodal graphs of diameter~$3$, is now the only class remaining where no general bounds (other than the trivial ones) on the metric dimension are known.  Finding such bounds is an interesting open problem!

\proof[Acknowledgements]
The author acknowledges financial support from an NSERC Discovery Grant and a Memorial University of Newfoundland startup grant.  He would like to thank Daniele Bartoli and Tam\'as H\'eger for a useful discussion.

\end{document}